\documentclass[12pt]{article}

\newcommand{\mathfrak}{\mathcal}
\newcommand{\mathbb}{\mathbf}

\usepackage{graphicx}

\title{Differential Algebra Structures \protect\\ on Families of Trees}
\author{Robert L. Grossman \and Richard G. Larson}
\date{June 6, 2004}

\newcommand{\D}{\mathfrak{D}}
\newcommand{\B}{\mathfrak{B}}
\newcommand{\fS}{\mathfrak{S}}

\newcommand{\T}{\mathcal{T}}

\newcommand{\cI}{\mathcal{I}}

\newcommand{\LT}[1]{\mathcal{T}({#1})}
\newcommand{\LOT}{\mathcal{LOT}}
\newcommand{\kLT}[1]{k\{\LT{#1}\}}
\newcommand{\Der}{\mathop{\mathrm{Der}}}
\newcommand{\DerR}{\mathop{\mathrm{Der}^R}}
\newcommand{\Diff}{\mathop{\mathrm{Diff}}}
\newcommand{\DiffR}{\mathop{\mathrm{Diff}^R}}
\newcommand{\Diffk}[1]{\mathop{\mathrm{Diff}_k}(\mathfrak{#1})}
\newcommand{\Ker}{\mathop{\mathrm{Ker}}}
\newcommand{\Imm}{\mathop{\mathrm{Im}}}
\newcommand{\eqdef}{\stackrel{\mathrm{def}}{=}}
\newcommand{\mapdef}[1]{\stackrel{#1}{\longrightarrow}}
\newcommand{\ptl}[1]{\partial/\partial{#1}}
\newcommand{\tensor}[1]{\mathbin{{\otimes_{#1}}}}
\newcommand{\mtimes}[1]{\mathbin{{\times_{#1}}}}
\newcommand{\btensor}{\tensor{\mathrm{r}}}
\newcommand{\smsh}[1]{\mathbin{{{\#}_{#1}}}}
\newcommand{\RR}{$(R,R)$}           
\newcommand{\pf}{\medskip\noindent{\sc Proof: }}
\newcommand{\pfend}{\medbreak}
\newtheorem{lemma}{Lemma}[section]
\newtheorem{dfn}[lemma]{Definition}
\newtheorem{prop}[lemma]{Proposition}
\newtheorem{thm}[lemma]{Theorem}

\newcommand{\Conn}[2]{\nabla_{#1}{#2}}
\newcommand{\uu}{\mathrm{u}}
\newcommand{\vv}{\mathrm{v}}
\newcommand{\ttt}{\mathrm{t}}
\newcommand{\EndR}{\mathrm{End}(R)}

   {\begin{list}{}{\setlength{\rightmargin}{\leftmargin}}}%
   {\end{list}}

\newenvironment{xmpl}%
 {\begin{trivlist}\refstepcounter{lemma}%
 \item[\hskip\labelsep{\bf Example\ \thelemma}]}{\end{trivlist}}

\newenvironment{const}%
 {\begin{trivlist}\refstepcounter{lemma}%
 \item[\hskip\labelsep{\bf Construction\ \thelemma}]}{\end{trivlist}}

\newcounter{atmz}

\newenvironment{atomize}%
   {\begin{list}{{\rm \alph{atmz})}}%
   {\usecounter{atmz}}}%
   {\end{list}}

\begin{document}

\maketitle

\begin{abstract}
It is known that the vector space spanned by labeled rooted trees forms
a Hopf algebra. Let
$k$ be a field and let $R$ be a commutative $k$-algebra.  Let $H$
denote the Hopf algebra of rooted trees labeled using derivations
$D\in \Der(R)$.  In this paper, we introduce a construction which
gives $R$ a $H$-module algebra structure and show this induces a
differential algebra structure of $H$ acting on $R$.  The work
here extends the notion of a $R/k$-bialgebra introduced by
Nichols and Weisfeiler.

\end{abstract}

\section{Introduction}

Let $k$ be a field, $R$ be a commutative $k$-algebra, and $\Der(R)$
the Lie algebra of derivations of $R$.  It is known that the vector
space spanned by labeled rooted trees forms a Hopf
algebra~\cite{GLtrees}.  Let $H$ denote the Hopf algebra of rooted
trees whose non root nodes are labeled using derivations $D\in
\Der(R)$~\cite{GLtrees}.  For such a Hopf algebra, we introduce a
class of $H$-module algebras which we call Leibnitz, and give a
construction which yields a variety of different Leibnitz $H$-module
algebras (Theorem~\ref{LeibMeasProp}).

We also show how Leibnitz $H$-module algebras are related to
Nichols and Weisfeiler's $R/k$-bialgebras~\cite{NichWeis},
which arise in Hopf-algebra approaches to differential algebra
(Theorem~\ref{LastThm}).
In Section~\ref{QuotSect} we also give a method for describing
quotients of Leibnitz $H$-module algebras (Theorem~\ref{QuotI0Thm}).

Hopf algebras can be used to simplify computations of derivations
\cite{GL:jsc-trees}.  In the same way, Leibnitz $H$-module algebras
can be used to simplify the symbolic computation of derivations acting
on polynomials and other algebras of functions.  This was used in
\cite{CG:manifolds} and \cite{CGL:geom-stable} to derive geometrically
stable numerical integration algorithms, although the results are
presented differently.

\section{Bialgebras and trees}\label{TreeSect}

This section reviews some material from~\cite{GLtrees} on the Hopf
algebra structure of trees.
For background material on Hopf algebras see~\cite{moss}.
Throughout this paper $k$ is a field.

Let $R$ be a commutative $k$-algebra.
Let $\D$ be a vector space over $k$, and let $\fS$ be a subset of $\D$.
Let $\LT{\fS}$ denote the set of
ordered trees in which each node other than the root is labeled with an
element of $\fS$.
Let $\kLT{S}$ be the vector space over $k$ with basis the
elements of $\LT{S}$.

When we speak of a ``subtree of a tree'' we also include the
(unlabeled) node
in the tree to which the subtree is attached as the root of the subtree.
When we refer to the ``children of a  node $v$'' we will sometimes mean
the nodes which are attached to $v$ as immediate descendants, and
sometimes mean the full subtrees which are attached to $v$.  Which sense
is meant will be clear from context.

We have a grading on $\kLT{\fS}$:
$\kLT{\fS}_r$ is spanned by the trees with $r+1$ nodes.
There is a bialgebra structure on $\kLT{\fS}$ given in~\cite{GLtrees}
which we summarize here.

\begin{figure}[hp]
\begin{center}
\includegraphics[scale=0.5]{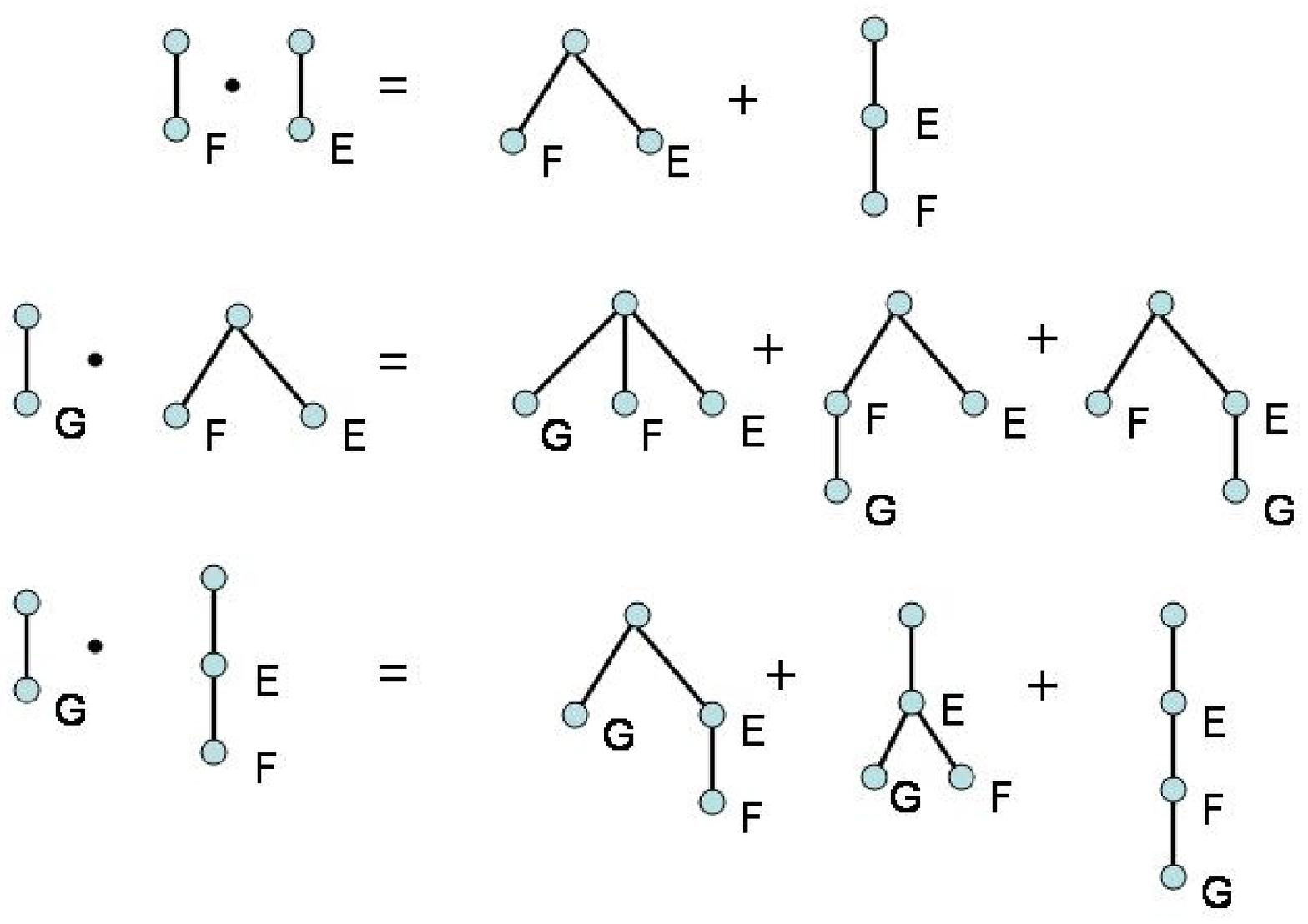}
\end{center}
\caption{Some examples of multiplying labeled trees.}
\label{fig:mul}
\end{figure}

Multiplication in $\kLT{\fS}$ is defined as follows.
Let $T_1$ and $T_2$ be rooted trees.
Remove the root from the tree $T_1$ to form a multiset $\mathcal{F}$ of
rooted trees.
Let $d$ be a function from $\mathcal{F}$ to the set of nodes of $T_2$.
Let $T_d$ be the rooted tree formed by adding an edge to link
the root of $T'\in\mathcal{F}$ to the node $d(T')$ of $T_2$.
The order on the children of a node of $T_d$ is given by declaring
that the order of the children of a node of a tree in $\mathcal{F}$ is
preserved,
that the order of the children of a node of $T_2$ is preserved,
that if two roots of trees in $\mathcal{F}$ are linked to the same node
of  $T_2$,
they are given the order they had as children of the root of $T_1$, and
that the root of $T'$ precedes every child of the node $d(T')$
to which it is linked.
Now
\[
T_1\cdot T_2 = \sum_d T_d,
\]
where the sum ranges over all trees $T_d$ formed as described above.
This product is associative, and the tree
with only one node is a multiplicative unit.
See~\cite{GLtrees} for details and Figure~\ref{fig:mul} for some examples.
(Note that in~\cite{GLtrees} multiplication is defined so that the root
of $T'\in\mathcal{F}$ follows every child of the node $d(T')$; this
will not change the application of the results we use
from there.)

\begin{figure}[hp]
\begin{center}
\includegraphics[scale=0.5]{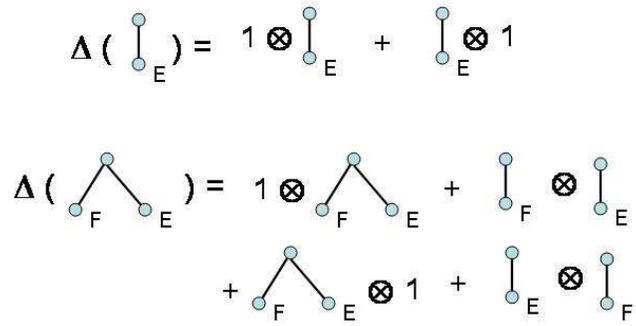}
\end{center}
\caption{Some examples of co-multiplying labeled trees.}
\label{fig:comul}
\end{figure}

The coproduct in $\kLT{\fS}$ is defined as follows.
If $T$ is a labeled ordered rooted tree, define
\[
\Delta(T)
= \sum_\mathcal{X} T_\mathcal{X}\otimes
                   T_{\mathcal{F}\backslash\mathcal{X}},
\]
where $\mathcal{X}$ ranges over all sub-multisets of the ordered
multiset $\mathcal{F}$ described above.
If $\mathcal{Y}$ is a sub-multiset of $\mathcal{F}$,
the labeled ordered tree $T_\mathcal{Y}$ is formed by adding edges to
link the roots of
the trees in $\mathcal{Y}$ to a new root, preserving their original
order and labels.
In particular, $T_\mathcal{F}=T$ and $T_\emptyset=1$ is the unit.
The counit $\epsilon : \kLT{\fS}\rightarrow k$ is defined as follows.
\[
\epsilon(T) = \left\{
\begin{array}{ll}
1 & \mbox{if $T = 1$}\\
0 & \mbox{otherwise.}
\end{array}
\right.
\]
In this construction, each non root node retains its original label;
recall that the root is not labeled.
See~\cite{GLtrees} for more details and Figure~\ref{fig:comul} for some
examples.

Since $\kLT{\fS}$ is a graded bialgebra with $\dim\kLT{\fS}_0=1$, it is
a Hopf algebra.


We summarize this discussion in the following

\begin{prop}\label{kLTSbialg}
Let $k$ be a field and let $\fS$ be a vector space over $k$.
Then $\kLT{\fS}$ is a cocommutative graded connected Hopf algebra.
\end{prop}

\section{$H$-module algebras}\label{MeasSect}

Let $R$ be a commutative $k$-algebra, and let $H$ be a $k$-bialgebra.
Recall the following definition~\cite[4.1.1]{Montgomery}:

\begin{dfn}\label{HModuleAlg}
The algebra $R$ is a {\em left $H$-module algebra\/} if $R$ is a left
$H$-module for which
\[
h\cdot(rs)=\sum_{(h)} (h_{(1)}\cdot r)(h_{(2)}\cdot s),
\]
where $h\in H$, $\Delta(h)=\sum_{(h)} h_{(1)}\otimes h_{(2)}$,
and $r$, $s\in R$.
\end{dfn}

\begin{xmpl}\label{ex3}
Let $R=k[X_1$, \ldots, $X_N]$.
Then the Lie algebra $\D$ of derivations of $R$ has $\{\ptl{X_1}$,
\ldots, $\ptl{X_N}\}$ as an $R$-basis.
Let $H$ be the bialgebra $\kLT{\D}$ defined in Section~\ref{TreeSect}.
We define an $H$-module algebra structure on $R$ as follows.
Let $T\in\kLT{\D}$ be a labeled tree.
Number the root of $T$ with~0, and number the other nodes
1, \ldots, $m$.
Let node $i$, $i>0$, be labeled with
\begin{equation}\label{ex3eq}
E_i = \sum_{\mu_i=1}^N r_{i\mu_i}\frac{\partial}{\partial X_{\mu_i}},
\end{equation}
where $r_{i\mu_i}\in R$.
Suppose that node $i$, $i\ge 0$, has children $j_1$, \ldots, $j_k$.
Define
\[
R(i;\mu_{j_k}\ldots\mu_{j_1}) = \left\{\begin{array}{ll}
  {\displaystyle\frac{\partial}{\partial X_{\mu_{j_k}}}}\cdots
    {\displaystyle\frac{\partial}{\partial X_{\mu_{j_1}}}s}
    & \mbox{if $i=0$,} \\*[6pt]
  {\displaystyle\frac{\partial}{\partial X_{\mu_{j_k}}}}\cdots
    {\displaystyle\frac{\partial}{\partial X_{\mu_{j_1}}}r_{i\mu_i}}
      & \mbox{otherwise.}
\end{array}\right.
\]
We will usually abbreviate $R(i;\mu_{j_k}\ldots\mu_{j_1})$ by
$R(i)$.
Define
\[
T\cdot s = \sum_{\mu_1,\ldots,\mu_m=1}^N R(m)\cdots R(1)R(0),
\]
for $s\in R$.
It can be shown~\cite{GL:am-solv}[Prop.~2] that this makes $R$ into a
left $H$-module algebra.
We will prove the existence of more complex $H$-module algebra
structures in Section~\ref{QuotSect}.
\end{xmpl}

We will generalize Example~\ref{ex3} in Proposition~\ref{mainthm}.
Here is a specific case of the construction in
Example~\ref{ex3}.

\begin{figure}[hp]
\begin{center}
\includegraphics[scale=0.5]{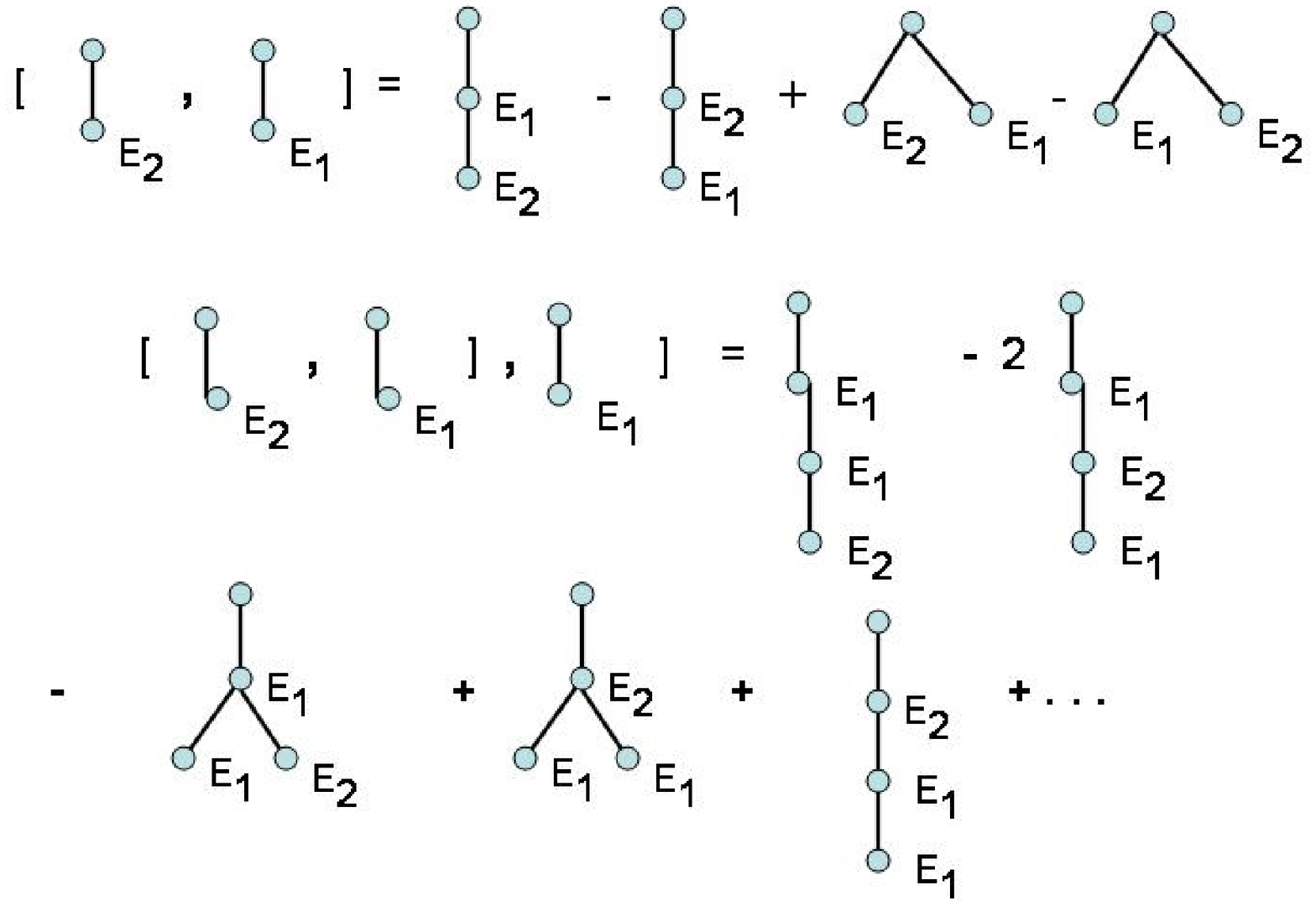}
\end{center}
\caption{This figure illustrates some of the trees which
arise in Example~\ref{exHall}.  The Lie bracket $[E_2, E_1], E_1]$
only shows some of the trees which arise in the expansion.
}
\label{fig:hall-ex1}
\end{figure}

\begin{xmpl}\label{exHall}
Consider the following two vector fields on
$\bf{R}^8$ introduced in~\cite{Grossman:1991}:
\begin{eqnarray*}
E_1 & = & {{\partial}\over {\partial x_1}} \\
E_2 & = & {{\partial}\over {\partial x_2}}
- x_1 {{\partial}\over {\partial x_3}}
+ {1\over 2} x_1^2 {{\partial}\over {\partial x_4}}
+  x_1 x_2 {{\partial}\over {\partial x_5}} \\
&& \quad {} - {1\over 6}x_1^3 {{\partial}\over {\partial x_6}}
- {1\over 2}x_1^2 x_2 {{\partial}\over {\partial x_7}}
- {1\over 2}x_1 x_2^2 {{\partial}\over {\partial x_8}}.
\end{eqnarray*}
Then it is simple to check, for example, that
\begin{eqnarray*}
[E_2, E_1] & = &
{{\partial}\over {\partial x_3}}
-x_1 {{\partial}\over {\partial x_4}}
-x_2 {{\partial}\over {\partial x_5}}
+ {1\over 2} x_1^2 {{\partial}\over {\partial x_6}}
+ x_1 x_2 {{\partial}\over {\partial x_7}}
+ {1\over 2} x_2^2 {{\partial}\over {\partial x_8}}, \\
{}[[E_2, E_1], E_1]  & = &
{{\partial}\over {\partial x_4}}
- x_1  {{\partial}\over {\partial x_6}}
- x_2  {{\partial}\over {\partial x_7}}, \\
{}[[E_2, E_1], E_2]  & = &
{{\partial}\over {\partial x_5}}
- x_1  {{\partial}\over {\partial x_7}}
- x_2  {{\partial}\over {\partial x_8}}.
\end{eqnarray*}
See Figure~\ref{fig:hall-ex1} for the corresponding trees.  Note also
that:
\[
[E_2, E_1](x_3) = 1, \quad
[[E_2, E_1], E_1](x_4) = 1, \quad
\mbox{ and } \quad [[E_2, E_1], E_2](x_5) = 1.
\]
We could continue this example by checking that the actions of the
differential operators is the same as the actions of the trees.
\end{xmpl}

We will use the following definition in the sequel.

\begin{dfn}\label{LilTreeDef}
Let $\D$ be a vector space.
\begin{atomize}
\item
Let $E\in \D$.
Denote by $\vv(E)$ the labeled ordered tree with two nodes:
the root, and a single child which is labeled with $E$.
\item
Let $T_1$, \ldots, $T_k$ be labeled ordered trees, and let $E\in\D$.
Denote by $\uu(E;T_1,\ldots,T_k)$ the labeled ordered tree whose root
has one child, labeled with $E$, with which the roots of the subtrees
$T_1$, \ldots, $T_k$ are identified.
The ordering on the children of the node labeled with $E$ in
$\uu(E;T_1,\ldots,T_k)$ is given by specifying that the children of
the root of $T_1$ precede the children of the root of $T_2$, \dots,
which precede the children of the root of $T_k$.
\item
Let $T_1$, \ldots, $T_k$ be labeled ordered trees.
Denote by $\ttt(T_1,\ldots,T_k)$ the labeled ordered tree formed by
identifying the roots of the trees $T_1$, \ldots, $T_k$.
The ordering on the children of the root $\ttt(T_1,\ldots,T_k)$ is given
by specifying that
the children of the root of $T_1$ precede the children of the root of
$T_2$,
\dots, which precede the children of the root of $T_k$.
\end{atomize}
\end{dfn}

\begin{figure}[hp]
\begin{center}
\includegraphics[scale=0.5]{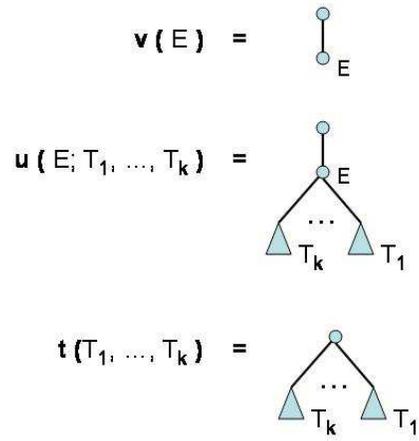}
\end{center}
\caption{This figure illustrates the maps $\vv(E)$,
$\uu(E;T_1,\ldots,T_k)$ and $\ttt(T_1,\ldots,T_k)$ which are
used to define an action of the algebra of labeled trees
$\kLT{\D}$ on $R$.
}
\label{fig:uvt}
\end{figure}

\begin{figure}[hp]
\begin{center}
\includegraphics[scale=0.5]{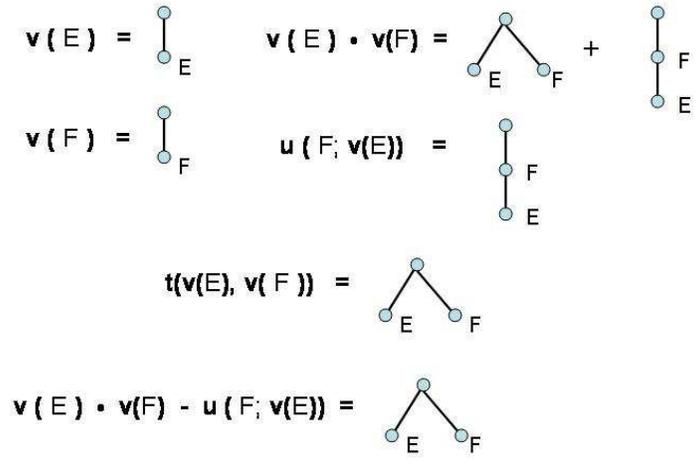}
\end{center}
\caption{This figure illustrates the computation
in Lemma~\ref{PairEquation}.
}
\label{fig:uvteq}
\end{figure}

This definition is illustrated in Figure~\ref{fig:uvt}.
Note that $\uu(E;T_1,\ldots,T_k)=\uu(E;\ttt(T_1,\ldots,T_k))$.
We will need the following lemma.

\begin{lemma}\label{PairEquation}
Let $E$, $F\in\D$.  Then
\[
\ttt(\vv(E),\vv(F))= \vv(E)\cdot\vv(F) - \uu(F;\vv(E)).
\]
\end{lemma}

\pf The proof of the lemma follows immediately from the
definition of multiplication for trees.

This Lemma is illustrated in Figure~\ref{fig:uvteq}.
\pfend

\begin{dfn}\label{conndef}
Let $R$ be a commutative $k$-algebra,
and let $\D$ be a Lie algebra of derivations of $R$.
A map $\D\times\D\rightarrow
\D$ sending $(E,F)\in\D\times\D$ to $\Conn{E}{F}\in\D$
satisfying
\begin{atomize}
\item
$\Conn{E_1+E_2}{F} = \Conn{E_1}{F}+\Conn{E_2}{F}$
\item
$\Conn{E}{(F_1+F_2)} = \Conn{E}{F_1}+\Conn{E}{F_2}$
\item
$\Conn{f\cdot E}{F} = f\cdot\Conn{E}{F}$
\item
$\Conn{E}{(f\cdot F)} = f\cdot\Conn{E}{F}+E(f)F$
\end{atomize}
where $E$, $F\in\D$, $f\in R$, is called a
{\em connection.}
\end{dfn}

(For example, $M$ could be a $C^\infty$-manifold,
$R$ be the algebra of $C^\infty$ functions on $M$,
$\D$ could be the Lie algebra of vector fields on $M$, and
$\Conn{E}{F}$ could be the Koszul connection (see~\cite{SpivakI},
Chapter~5, and~\cite{SpivakII}, Chapters~5 and~6).)

\begin{figure}[hp]
\begin{center}
\includegraphics[scale=0.5]{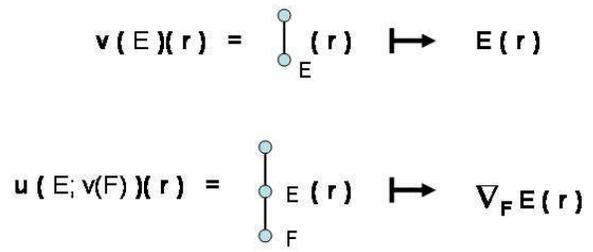}
\end{center}
\caption{An illustration of the action defined by Construction~\ref{ex4}.}
\label{fig:action1}
\end{figure}

\begin{figure}[hp]
\begin{center}
\includegraphics[scale=0.5]{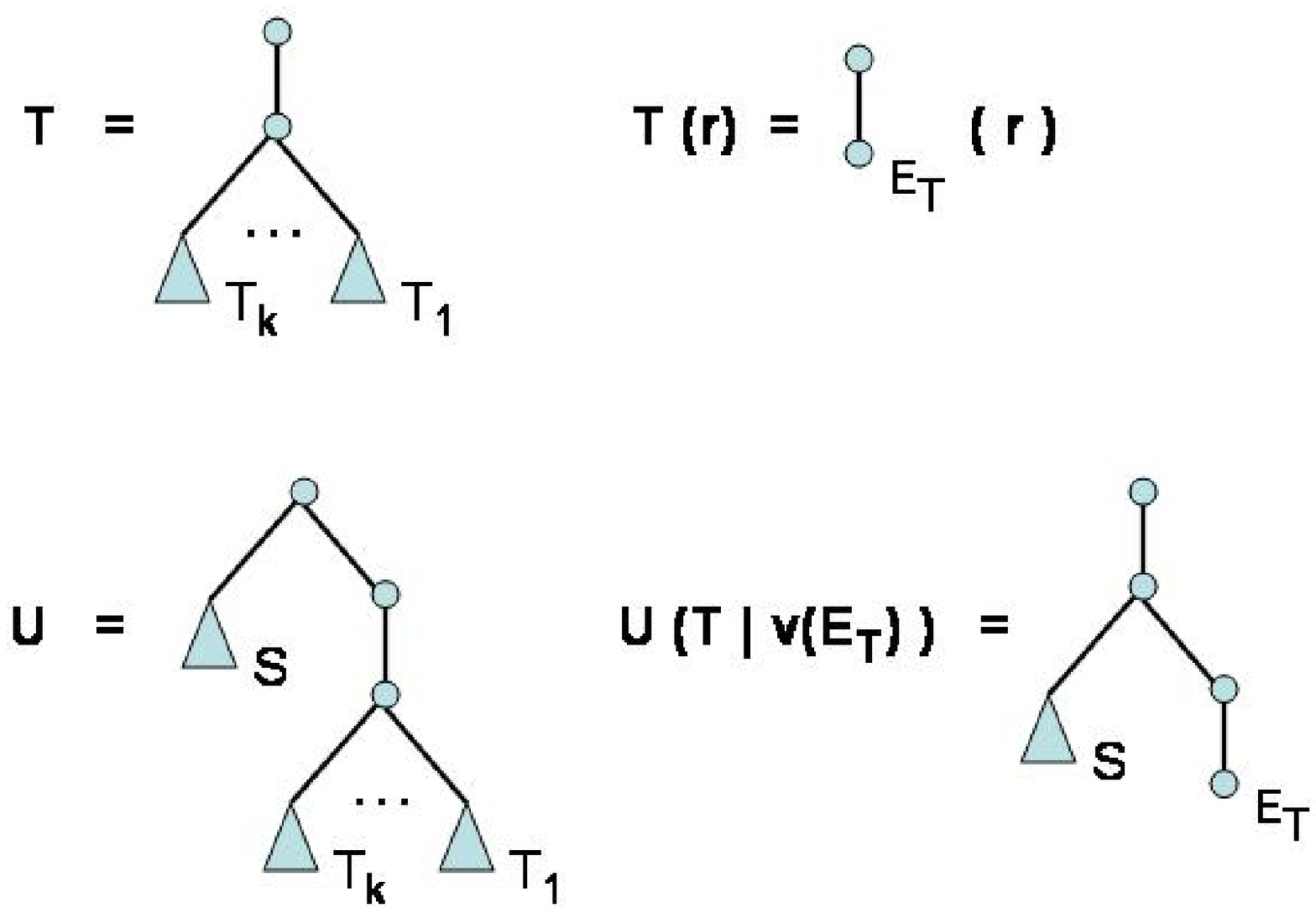}
\end{center}
\caption{Another illustration of the action defined by Construction~\ref{ex4}.}
\label{fig:action2}
\end{figure}

\begin{const}\label{ex4}
We use the action of $\D$ on $R$ and a connection
$\D\times\D\rightarrow\D$ to construct an action
of the algebra of labeled trees $\kLT{\D}$ on $R$.
See Figures~\ref{fig:action1} and \ref{fig:action2}.
We give an inductive description of the action.
The description of this construction is intended to allow an inductive
proof of Theorem~\ref{mainthm}.
We make the following assumptions about the action.
\begin{atomize}
\item\label{OpReplace} The tree $\vv(E)$ acts as $E$.
\item\label{ConnReplace} The tree $\uu(E;\vv(F))$ acts as $\Conn{F}{E}$.
\item\label{CohReplace}
Suppose that $T$ is a labeled ordered tree whose root has a single child
and which acts on $R$ as the differential operator $E_T$.
Suppose further that
$U$ is a labeled ordered tree which contains $T$ as a proper subtree.
Denote by $U(T|\vv(E_T))$ the labeled ordered tree resulting from
replacing the subtree $T$ with the tree $\vv(E_T)$.
In this construction, we require that $U$ acts like $U(T|\vv(E_T))$.

This assumption says that a subtree whose root has one child can be
replaced by a tree which has one non root node which is labeled with a
derivation whose action is that of the original subtree.
\end{atomize}
We make use of these assumptions as follows.
\begin{atomize}
\setcounter{atmz}{3} 
\item\label{PairReplace} If $E$, $F\in\D$, by
Lemma~\ref{PairEquation} the tree
$\ttt(\vv(E),\vv(F))$ acts as $\vv(E)\cdot\vv(F)-\uu(F;\vv(E))$,
whose action we know by~(\ref{OpReplace}) and~(\ref{ConnReplace}).
\item\label{PrimReplace} We define the action of a tree whose root
has only one
child by induction on the number of children of the child of the root:
\begin{itemize}
\item If the child of the root has one child, by induction on the number
of nodes we know how that child acts, since
application of~(\ref{CohReplace}) and~(\ref{ConnReplace}) allows
us to determine the action of the tree on $R$.
\item Suppose that $T=\uu(E;T_1,\ldots,T_{n+1})$, where each $T_i$ is a
tree whose root has only one child. Then
\[
T_1\cdot\uu(E;T_2,\ldots,T_{n+1})=\ttt(T_1,\uu(E;T_2,\ldots,T_{n+1}))+
T+\sum_j U_j,
\]
where the $U_j$ are trees in which $T_1$ has been linked to various
nodes (other than the root) in the trees $T_j$ in
$\ttt(E;T_2,\ldots,T_{n+1})$.
By induction we know the action of $T_1$, of
$\uu(E;T_2,\ldots,T_{n+1})$, and of the $U_j$ on R.
By~(\ref{CohReplace}) and~(\ref{PairReplace}) we know the action of
\newline
$\ttt(T_1,\uu(E;T_2,\ldots,T_{n+1}))$.
Therefore we know the action of $T$.
\end{itemize}
This gives the action of a tree whose root has only one child.
\item~\label{GenlReplace} We now determine the action of a general tree
by induction on the number of children of the root:
\begin{itemize}
\item The case where the root has one child follows from~(\ref{PrimReplace}).
\item Suppose that $T=\ttt(T_1,\ldots,T_{n+1})$, where each $T_i$ is a
tree whose root has one child.
Now
\[
T_1\cdot\ttt(T_2,\ldots,T_{n+1})= T+\sum_j V_j,
\]
where the $V_j$ are trees whose roots have $n$ children.
We know the action of $T_1$ by~(\ref{PrimReplace}), and of
$\ttt(T_2,\ldots,T_{n+1})$ and $V_j$ by induction.
Therefore we know the action of $T$.
\end{itemize}
\end{atomize}
\par\noindent
Note that this construction includes Example~\ref{ex3}, upon letting
\[
\Conn{E_i}{E_j}= \sum_{\mu_i,\mu_j} r_{i\mu_i}
\frac{\partial r_{j\mu_j}}{\partial X_{\mu_i}}\frac{\partial}
{\partial X_{\mu_j}}
\]
where the $E_i$ are the differential operators given in~(\ref{ex3eq}).
This follows from Definition~\ref{conndef} and the fact
that $\Conn{\frac{\partial}{\partial X_i}}
{\displaystyle\frac{\partial}{\partial X_j}} = 0$.
\end{const}

Note that Construction~\ref{ex4} gives the basic description
of the action of a bialgebra of ordered trees whose non root nodes are
labeled with derivations of $R$  on the commutative algebra $R$.

\begin{dfn}
If $T_1$, $T_2\in\kLT{\D}$ act identically on $R$
for $T_1$, $T_2\in\kLT{\D}$,
write $T_1\sim T_2$.
\end{dfn}

\begin{thm}\label{mainthm}
Let $R$ be a commutative $k$-algebra, and let $\Conn{E}{F}$ be a
connection on the Lie algebra $\D$ of derivations of $R$.
Then Construction~\ref{ex4} gives a $\kLT{\D}$-module structure on $R$.
This module structure induces a map $\psi:\kLT{\D}\longrightarrow\EndR$.
The following conditions are satisifed:
\begin{enumerate}
\item $h\cdot r=\psi(h)\cdot r$ for all $h\in\kLT{\D}$, $r\in R$;
\item the $\kLT{\D}$-module structure on $R$ is a
$\kLT{\D}$-module algebra structure;
\item $\Imm\psi\subseteq\Diff(R)$.
\end{enumerate}
\end{thm}

\pf
\begin{enumerate}

\item\label{No2} Assumption~(\ref{OpReplace})
in Construction~\ref{ex4} guarantees that the actions of $h$ and
$\psi(h)$ are the same if $h$ is a tree with one non root node.

Assumption~(\ref{ConnReplace}) describes how a tree with two nodes whose
root has one child acts.

Condition~(\ref{CohReplace}) says that subtrees which act as derivations
can be replaced by trees with one non root node which act as the same
derivation.

Condition~(\ref{PairReplace}) says that the actions of $h$ and
$\psi(h)$ are the same if $h$ is the product of two trees which have
only one non root node.

Conditions~(\ref{PrimReplace}) and~(\ref{GenlReplace}) are used to
prove the result by induction on the number of children of the root of
any tree, and on the number of children of the child of the root in the
case that the root has only one child.

\item The fact that the action gives a $\kLT{\D}$-module algebra
structure follows from the definition of the coalgebra structure of
$\kLT{\D}$.

\item $\Diff(R)$ is generated by derivations.
Since Construction~\ref{ex4} gives the action of any tree as a sum of
products of derivations, it follows that $\psi(h)\in\Diff(R)$.

\end{enumerate}

\pfend

\begin{dfn}\label{LeibMeasDef}
Let $T$ be a labeled ordered tree.
Suppose that a non root node $i$ of $T$ is labeled with $rE$, where
$r\in R$ and $E\in\D$.
Denote the labeled ordered subtree whose (unlabeled) root is $i$ by
$T_i$.
Denote by $T(i,G,T')$ the tree identical to $T$, except that the node
$i$ of $T(i,G,T')$ is labeled with $G\in\D$, and the labeled
ordered subtree rooted at node $i$ is $T'\in\kLT{\D}$.
(Note that $T(i,rE,T_i)=T$.)
Extending this notation allows us to replace $T_i$ with a linear
combination of trees.

Let $R$ be a commutative $k$-algebra, let $\D=\Der(R)$,
and suppose that $R$ is a Hopf module algebra over the $k$-bialgebra
$\kLT{\D}$.
The $\kLT{\D}$-module algebra structure of $R$ is called {\em
Leibnitz\/} if
\[
T\cdot s=\sum_{(T_i)}({T_i}_{(1)}\cdot r)(T(i,E,{T_i}_{(2)})\cdot s)
\]
for all trees $T$, for all non root nodes $i$ of $T$, and for all
factorizations
$rE$, where $r\in R$, $E\in\D$, of the label of node $i$.
(Note the coproduct in the formula above is over the tree $T_i$,
not over the tree $T$.)
\end{dfn}

\begin{figure}[hp]
\begin{center}
\includegraphics[scale=0.5]{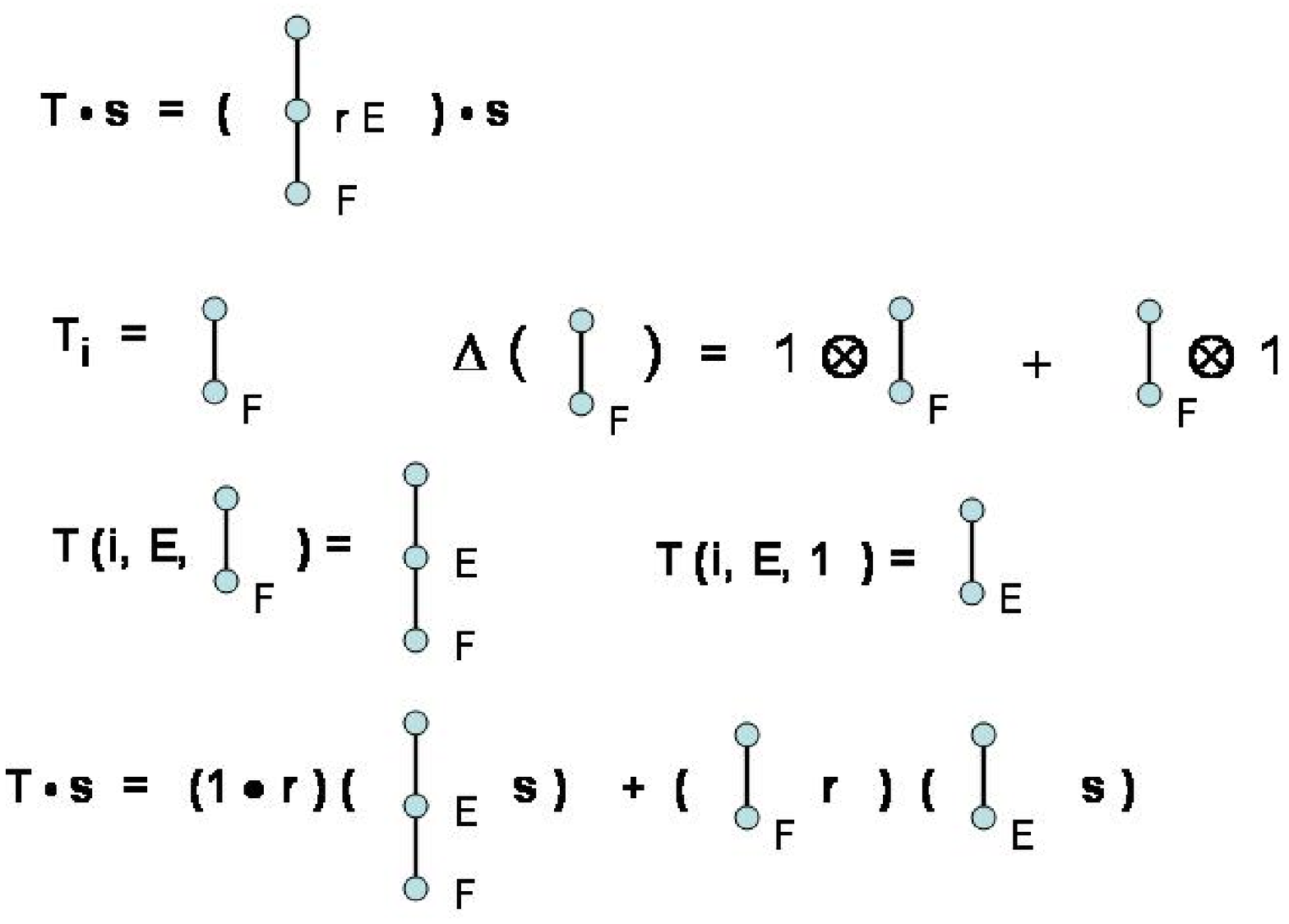}
\end{center}
\caption{This figure
illustrates the defining formula for
a Leibnitz $H$-module algebra structure for a simple
tree acting on a function $s$.  Note that in this simple case,
the formula defining a Leibnitz structure generalizes the following
standard formula for connections:
$(\Conn{F}{(rE)})(s) = r ((\Conn{F}{E})s) + F(r) E(s).$
}
\label{fig:leibnitz-ex1}
\end{figure}

Basically the $\kLT{\D}$-module algebra structure is Leibnitz if the
action of a labeled tree is consistent in that subtrees act consistently
with how they act as seperate trees.  See Figure~\ref{fig:leibnitz-ex1}
for a simple example.

\begin{thm}\label{LeibMeasProp}
Let $R$ be a commutative $k$-algebra, and $\D=\Der(R)$,
and let $\nabla:\D\times\D\rightarrow\D$ be a connection.
Then the $\kLT{\D}$-module algebra structure on $R$ given in
Construction~\ref{ex4} is Leibnitz.
\end{thm}

\pf
Let $T\in \kLT{\D}$, and suppose that the subtree rooted at node~$i$ is
$T_i=\ttt(T_{i1},\ldots,T_{ik})$, where each $T_{ij}$ is a tree whose
root has only one child.
By assumption~(\ref{CohReplace}) of Construction~\ref{ex4}
the subtree $T_{ij}$ can be replaced by $U_j=\vv(F_{j})$ where
$T_{ij}$ acts as $F_{j}$, so that
$T_i$ can be replaced by $U=\ttt(U_1,\ldots,U_k)$.

We show by induction on $k$ that
\[
T\cdot s = \sum_{(U)} (U_{(1)}\cdot r) (T(i,E,U_{(2)})\cdot s).
\]

If $k=1$ then let $V$ be
the subtree rooted at parent node of~$i$ which has
node~$i$ as its one child, which is labeled
with $rE$ and which has one child node labeled with $F_1$, that is, the
tree $\uu(rE;\vv(F_{1})$.
By assumption~(\ref{ConnReplace}) of Construction~\ref{ex4},
the subtree $V$ acts as
\[
\Conn{F_1}{(rE)} = r\Conn{F_1}{E}+F_1(r)E
\]
so that
\[
T\cdot s =\sum_{(U)}(U_{(1)}\cdot r)(T(i,E,U_{(2)})\cdot s)
\]
in this case.

Now consider the case $k>1$.
Define $G_\ell=\Conn{F_1}{F_\ell}$.
Here
\begin{eqnarray}
\lefteqn{\uu(rE;\vv(F_1),\ldots,\vv(F_k)) } \nonumber \\
& \sim & \vv(F_1)\cdot\uu(rE;\vv(F_2),\ldots,\vv(F_k)) \nonumber \\
& & \qquad {}-\ttt(\vv(F_1),\uu(rE;\vv(F_2),\ldots,\vv(F_k)) \nonumber \\
& & \qquad
{}-\sum_{\ell=2}^k\uu(rE;\vv(F_2),\ldots,\vv(G_\ell),\vv(F_k))
\nonumber \\
& \sim &
\vv(F_1)\cdot\left( \sum_{(V)}(V_{(1)}\cdot r)
 (\uu(E;V_{(2)})\right) \label{leib1} \\
& & \qquad {}-\sum_{(V)}(V_{(1)}\cdot r)(\ttt(\vv(F_1),\uu(E;V_{(2)})
\nonumber \\
& & \qquad {}-\sum_{\ell,(V)}(V_{\ell(1)}\cdot
r)(\uu(E;V_{\ell(2)}) \nonumber \\
& \sim &
\sum_{(U)}(U_{(1)}\cdot r)(\uu(E;U_{(2)}) \nonumber \\
& & \qquad {}+ \sum_{(V)}(V_{(1)}\cdot r)(\ttt(\vv(F_1),\uu(E;V_{(2)})
\nonumber \\
& & \qquad {}+ \sum_{\ell,(V)}(V_{\ell(1)}\cdot
r)(\uu(E;V_{\ell(2)}) \nonumber \\
& & \qquad {}-\sum_{(V)}(V_{(1)}\cdot r)(\ttt(\vv(F_1),\uu(E;V_{(2)})
\nonumber \\
& & \qquad {}-\sum_{\ell,(V)}(V_{\ell(1)}\cdot
r)(\uu(E;V_{\ell(2)}) \nonumber \\
& \sim &
\sum_{(U)}(U_{(1)}\cdot r)(\uu(E;U_{(2)}) \nonumber
\end{eqnarray}
(Note that we are using $V$ as a ``local variable'' in each sum, and
that the value varies, depending on context.)
The identity of the actions of the terms in expression~(\ref{leib1})
follows from the induction hypothesis.
Since $\uu(rE;\vv(F_1),\ldots,\vv(F_k))\sim
\sum_{(U)}(U_{(1)}\cdot r)(\uu(E;U_{(2)})$ the theorem now follows from
assumption~(\ref{CohReplace}) of Construction~\ref{ex4}.

\pfend

\section{$R/k$-bialgebras}\label{R/kSect}

The notions of {\em $R/k$-bialgebra\/} and {\em $R/k$-Hopf algebra\/} as
described in~\cite{NichKk} and~\cite{NichWeis} capture many of the
essential aspects of differential algebra.
We first review some of the material found there.

Let $R$ be a $k$-algebra.
A $R/k$-algebra is a $k$-algebra $B$ into which $R$ is embedded.
Note that this makes $B$ into a left and right $R$-module, and that
$(rb)s=r(bs)$ for all $r$, $s\in R$, $b\in B$,
and that also $(rb)c=r(bc)$, $(br)c=b(rc)$, and $(bc)r=b(cr)$
for all $b$, $c\in B$, $r\in R$.
When we refer to the $R$-module structure of $B$, we will understand the
left $R$-module structure.

We denote by $B\tensor{R}B$ the tensor product of $B$ with itself using
the left $R$-module structure of $B$.
That is, with $(rb)\otimes c=b\otimes(rc)$.
Note that in general $B\tensor{R}B$ is not an $R$-algebra.

The material we present here is related to the
$\times_R$-bialgebra construction given in~\cite{moss2}.
There a $\times_R$-bialgebra is defined in terms of maps between $B$ and
$B\mtimes{R}B$, but the space $B\mtimes{R}B$ in~\cite{moss2} is only a
$k$-subspace of $B\tensor{R}B$.

\begin{dfn}\label{R/kBialgDef}
A {\em $R/k$-bialgebra\/} is a $R/k$-algebra $B$ together with
$R$-module maps $\Delta : B\rightarrow B\tensor{R}B$ and
$\epsilon : B\rightarrow R$ satisfying
\begin{atomize}
\item
$B$ together with the maps $\Delta$ and $\epsilon$ is a coalgebra
over $R$.
\item
\label{PrCopr1} $\Delta(1)=1\otimes1$.
\item
\label{PrCopr2}
For all $b$, $c\in B$,
if $\Delta(b)=\sum\limits_{i}b_i\otimes b'_i$ and
$\Delta(c)=\sum\limits_{j}c_j\otimes c'_j$ are any representations of
$\Delta(b)$, $\Delta(c)\in B\tensor{R}B$, then
$\Delta(bc)=\sum\limits_{i,j}b_i c_j\otimes b'_i c'_j$.
\item
\label{UnCopr1} $\epsilon(1)=1$.
\item
\label{UnCopr2} $\epsilon(bc)=\epsilon(b\epsilon(c))$.
\end{atomize}
\end{dfn}

Note that condition~(\ref{PrCopr2}) of Definition~\ref{R/kBialgDef}
implies that
\[
\Delta(br)=\sum_{(b)}b_{(1)}r\otimes b_{(2)} =
\sum_{(b)}b_{(1)}\otimes b_{(2)}r
\]
for $b\in B$ and $r\in R$.
It can be shown
(see~\cite{NichKk} for details)
that conditions~(\ref{PrCopr1}) and~(\ref{PrCopr2}) are equivalent to
the assertion that the action of $B$ on $B\tensor{R}B$ defined by
$b\cdot(c\otimes d)\eqdef\Delta(b)(c\otimes d)$,
for $b$, $c$, $d\in H$,
gives $B\tensor{R}B$ a well-defined left $B$-module structure, and
that conditions~(\ref{UnCopr1}) and~(\ref{UnCopr2}) are equivalent to
the
assertion that the action of $B$ defined on $R$ by
$b\cdot r\eqdef\epsilon(br)$, for $b\in B$ and $r\in
R$, gives $R$ a well-defined left $B$-module structure.

If $B$ is a $R/k$-algebra, then $B$ is a \RR-bimodule via the left
and right actions of $R$ on $B$ induced by the embedding of $R$ in $B$.
Denote by $B\btensor B$ the tensor product of $B$ with itself
using this \RR-bimodule structure.
That is, $B\btensor B$ is an \RR-bimodule with
$r(b\btensor c)=(rb)\btensor c$, $(b\btensor c)r=b\btensor(cr)$,
and that $br\btensor c=b\btensor rc$.
The multiplication on $B$ induces a map $\mu:B\btensor B\rightarrow B$.

\begin{dfn}\label{R/kHopfAlgDef}
Let $B$ be a $R/k$-bialgebra.
An {\em antiproduct for $B$\/} is a $k$-linear map $E:B\rightarrow
B\btensor B$ satisfying
\begin{atomize}
\item
$E(rb)=rE(b)=E(b)r$ for all $r\in R$, $b\in B$;
\item
$\sum_{(b)} E(b_{(1)})b_{(2)} = b\btensor1$ for all $b\in B$;
\item
$(I\tensor{R}E)\circ\Delta(b)=(\Delta\btensor I)\circ E(b)$ for all
$b\in B$;
\item
$\mu\circ E(b)=\epsilon(b)1$ for all $b\in B$.
\end{atomize}
A $R/k$-bialgebra which has an antiproduct is called a {\em $R/k$-Hopf
algebra\/}.
\end{dfn}

We recall from~\cite[Proposition 8]{NichKk}

\begin{prop}\label{RHbialgProp}
Let $R$ be a $k$-algebra.
If $H$ is a cocommutative $k$-Hopf algebra
over which $R$ if a $H$-bimodule algebra,
then $R\smsh{k}H\eqdef R\tensor{k}H$
is a $R/k$-Hopf algebra,
with the $R$-coalgebra structure given by
\[
R\tensor{k}H\mapdef{{I_H}\otimes\Delta_H}
R\tensor{k}H\tensor{k}H \cong
(R\tensor{k}H)\tensor{R}(R\tensor{k}H)
\]
and
with multiplication given by
\[
(r\smsh{}h)(s\smsh{}k) \eqdef
\sum_{(h)}r(h_{(1)}\cdot s)\smsh{}h_{(2)}k
\]
where $r$, $s\in R$, $h$, $k\in H$,
and with antiproduct given by
\[
E(r\smsh{}h)=\sum_{(h)} (r\smsh{}h_{(1)})\btensor(1\smsh{}S(h_{(2)}))
\]
where $r\in R$, $h\in H$, and $S$ is the antipode of $H$.
\end{prop}

\begin{prop}\label{ConsistMorph}
Let $R$ be a $k$-algebra, let
$H$ and $\bar{H}$ be cocommutative $k$-Hopf algebras over which $R$ is a
module algebra,
and let $\varphi:H\rightarrow\bar{H}$ be a $k$-Hopf
algebra homomorphism which is consistent with the
$H$ and $\bar{H}$-module algebra structures on $R$.
Then the map
\[
R\smsh{k}H\rightarrow R\smsh{k}\bar{H}
\]
given by $r\smsh{}h\mapsto r\smsh{}\varphi(h)$ is a $R/k$-Hopf
algebra homomorphism.
\end{prop}

\pf Omitted.
\pfend

\begin{thm}\label{TreesRkThm}
Let $R$ be a commutative $k$-algebra, let $\D=\Der(R)$,
and let $\fS\subseteq\D$.
\begin{atomize}
\item
Let $\Diffk{S}$ denote the $k$-algebra of
higher order derivations generated by $\fS$, that is, the subalgebra
of $U(\D)$ generated by $\fS$.
Then the smash product $R\smsh{k}\Diffk{S}$ is a $R/k$-Hopf algebra.
The subalgebra of $\Diff(R)$ generated by $\fS$ is a homomorphic image
of $R\smsh{k}\Diffk{S}$.
\item
Let $\kLT{\fS}$ denote the Hopf
algebra of trees labeled with elements of $\fS$
defined in Section~\ref{TreeSect}.  Then
the smash product $R\smsh{k}\kLT{\fS}$ is a $R/k$-Hopf algebra.
\item\label{ItemC} The map
\[
R\smsh{k}\kLT{\fS}\rightarrow R\smsh{k}\Diffk{S}
\]
is a $R/k$-bialgebra homomorphism.
\end{atomize}
\end{thm}

\pf
The action of $\kLT{\D}$ on $R$ given in Construction~\ref{ex4} defines
a
Hopf algebra homomorphism $\kLT{\D}\rightarrow\Diff(R)$ which is
consistent with the Hopf module-algebra structures on $R$, and so
induces a
$R/k$-Hopf algebra homomorphism
$R\smsh{k}\kLT{\D}\rightarrow R\smsh{k}\Diff(R)$,
which allows us to use the differential algebra structure on
$R\smsh{k}\kLT{\D}$ to study the differential algebra structure on
$R\smsh{k}\Diff(R)$

$\Diffk{S}$ is a cocommutative $k$-Hopf algebra, so $R\smsh{k}\Diffk{S}$
is a $R/k$-Hopf algebra.

By Proposition~\ref{kLTSbialg}, $\kLT{\fS}$
is a cocommutative $k$-Hopf algebra.
Therefore, by Proposition~\ref{RHbialgProp}, $R\smsh{k}\kLT{\fS}$
is a $R/k$-Hopf algebra.

Part~(\ref{ItemC}) follows immediately from Proposition~\ref{ConsistMorph}.
\pfend

If $R$ is a $k$-bialgebra (for example if $R$ is the coordinate ring of
an affine algebraic group), then
the set $\DiffR(R)$ of right invariant differential operators is a
cocommutative $k$-Hopf algebra.
The following proposition follows immediately
from~\cite[Theorem 2.4.5]{MossBob}.

\begin{prop}\label{RtInvDiffProp}
Let $R$ be a $k$-bialgebra, and let $\DiffR(R)$ be the $k$-Hopf algebra
of right invariant differential operators.
Then
\[
\Diff(R)\cong R\smsh{k}\DiffR(R)
\]
is a $R/k$-Hopf algebra.
\end{prop}

If $R$ is a $k$-bialgebra, and if
$\B$ is a $k$-basis for $\DerR(R)$, the Lie algebra of right
invariant derivations of $R$ to itself, then we will see that
Proposition~\ref{RtInvDiffProp} allows us to use
the smash product $R\smsh{k}\kLT{\B}$ to do formal computations in
$\Diff(R)$.
(For example, such a basis always exists  if $R$ is the coordinate ring
of an affine algebraic group.)

We can use the $R/k$-bialgebra $R\smsh{k}\kLT{\fS}$ to do formal
computations involving elements of a subset $\fS\subseteq\Der(R)$.

Let $\fS\subseteq\Der(R)$, and let $F(\fS)$ be the free associative
algebra generated by $\fS$.
Recalling that the elements of $\fS$ are primitive, we get a Hopf
algebra structure on $F(\fS)$.
Since $F(\fS)$ is freely generated by $\fS$, we have maps
$F(\fS)\rightarrow\Diffk{S}$ mapping $E\in\fS$ to $E\in\Diffk{S}$,
and $F(\fS)\rightarrow\kLT{\fS}$ mapping $E\in\fS$ to
$\vv(E)\in\kLT{\fS}$.
These induce maps $p:R\smsh{k}F(\fS)\rightarrow R\smsh{k}\Diffk{S}$
and $i:R\smsh{k}F(\fS)\rightarrow R\smsh{k}\kLT{\fS}$.
If $R$ is an algebra for which there is a
connection on $\Der(R)$ (for example, if $R$ is the algebra of
$C^\infty$ functions on a Riemannian manifold), there is a map
$\varphi:R\smsh{k}\kLT{\fS}\rightarrow R\smsh{k}\Diffk{S}$ induced by
the map described in Construction~\ref{ex4}.
We have

\begin{thm}\label{LastThm}
Let $R$ be a commutative algebra for which there is a connection
on $\Der(R)$, and let $p$, $i$, and $\varphi$ be the maps described
above.
Then the diagram
\begin{center}
\begin{picture}(170,100)
\put(0,0){$R\smsh{k}F(\fS)$}
\put(55,80){$R\smsh{k}\kLT{\fS}$}
\put(130,0){$R\smsh{k}\Diffk{S}$}
%
\put(25,15){\vector(1,1){55}}
\put(48,45){$i$}
\put(100,70){\vector(1,-1){55}}
\put(130,45){$\varphi$}
\put(61,3){\vector(1,0){65}}
\put(90,9){$p$}
\end{picture}
\end{center}
commutes.
\end{thm}

This theorem allows us to do formal computations involving elements of
$\fS$ in the algebra $R\smsh{k}\kLT{\fS}$ rather than in
$R\smsh{k}F(\fS)$.

\section{Quotients of $R/k$-bialgebras}\label{QuotSect}

In this section we discuss certain quotients of the $R/k$-Hopf algebra
$R\smsh{k}\kLT{\D}$, where $\D=\Der(R)$.
The main result in this section is Theorem~\ref{QuotI0Thm}, which says
that a Leibnitz action of the $R/k$-Hopf algebra $R\smsh{k}\kLT{\D}$
can be computed from the action of $R\smsh{k}\kLT{\B}$ if $\B$ is an
$R$-basis of $\D$.

Let $B$ be a $R/k$-bialgebra.
A {\em $R/k$-biideal\/} is an ideal $I$ in the $R/k$-algebra $B$, such
that $\epsilon(I)=0$, and such that if $\pi:B\rightarrow B/I$ is the
projection of
$B$ onto $B/I$, we have $(\pi\tensor{R}\pi)\circ\Delta(I)=0$.
If $I$ is a $R/k$-biideal in $B$, then $B/I$ is a
$R/k$-bialgebra.
If $B$ is a $R/k$-Hopf algebra with antiproduct $E$, and if
$(\pi\btensor\pi)\circ E(I)=0$, then $I$ is called a {\em $R/k$-Hopf
ideal\/}, and $B/I$ is a $R/k$-Hopf algebra.

We will use the following definition in the sequel.

\begin{dfn}\label{I_0dfn}
Let $R$ be a Leibnitz $\kLT{\D}$-module algebra.
\newline
Let $\cI(\kLT{\D})$ be the $R$-linear span
of the elements of the form
\begin{equation}\label{I0Gen}
1\smsh{}T - \sum_{(T_i)}{T_i}_{(1)}\cdot r\smsh{}T(i,E,{T_i}_{(2)}),
\end{equation}
where $T\in\kLT{\D}$, with non root
node $i$ labeled with $rE$, with $r\in R$, $E\in\D$ (we include all
possible factorizations of the label of node $i$).
\end{dfn}

\begin{lemma}\label{LinSpanLemma}
The subspace $\cI(\kLT{\D})$ defined in Definition~\ref{I_0dfn} is a
two-sided ideal in $R\smsh{k}\kLT{\D}$.
\end{lemma}

\pf
Let $J$ be the ideal of $R\smsh{k}\kLT{\D}$ generated by
$\cI(\kLT{\D})$.

Let $Z$ be an element of the form~(\ref{I0Gen}), and let $T'$ be any
tree.
It follows from the definition of the product of trees that
$Z(1\smsh{}T')$
is an $R$-linear combination of elements of the form~(\ref{I0Gen}).

We now show that
\begin{eqnarray*}
(1\smsh{}T')Z & = &
1\smsh{}T'\cdot(1\smsh{}T - \sum_{(T_i)}{T_i}_{(1)}\cdot
      r\smsh{}T(i,E,{T_i}_{(2)})) \\
& = & 1\smsh{}T'\cdot T -
      \sum_{(T'),(T_i)}{T'}_{(1)}\cdot{T_i}_{(1)}\cdot
      r\smsh{}{T'}_{(2)}\cdot T(i,E,{T_i}_{(2)})
\end{eqnarray*}
is an $R$-linear combination of elements of the form~(\ref{I0Gen}).
Since $\kLT{\D}$ is generated as an algebra by trees whose root has one
child (see \cite{GLtrees}), it is sufficient to show this in the case
that the root of the tree $T'$ has one child.
Terms in the tree product $T'\cdot T$
pair in an obvious fashion with terms in
$\sum_{(T'),(T_i)}{T'}_{(1)}\cdot{T_i}_{(1)}\cdot
r\smsh{}{T'}_{(2)}\cdot T(i,E,{T_i}_{(2)})$.
Therefore $(1\smsh{}T')\cdot Z$ is again
an $R$-linear combination of elements of the form~(\ref{I0Gen}).
This completes the proof of the lemma.

\pfend

\begin{prop}\label{I0HopfIdlProp}
Let $R$ be a commutative $k$-algebra, let $\D=\Der(R)$,
and suppose that $R$ is a Leibnitz $\kLT{\D}$-module algebra.
\newline
Then the ideal $\cI(\kLT{\D}$ is a $R/k$-Hopf ideal.
\end{prop}
\pf

It is immediate that $\epsilon$ is zero on any element of the
form~(\ref{I0Gen}), since it is zero on any tree with more than one
node.

Let $\pi$ be the projection of $\kLT{\D}$ onto $\kLT{\D}/\cI(\kLT{\D}$.
To see that $(\pi\tensor{R}\pi)\circ\Delta$ is zero
on any element of the form~(\ref{I0Gen}),
write $\Delta(T)=\sum_j T'_j\otimes T''_j\in\kLT{\D}\tensor{k}\kLT{\D}$.
If node $i$ occurs in $T'_j$, then the corresponding term arising in the
coproduct applied to element~(\ref{I0Gen}) is
\begin{eqnarray*}
\lefteqn{1\smsh{}T'_j\tensor{R}1\smsh{}T''_j -
\sum_{(T_{i})}{T_{i}}_{(1)}\cdot r\smsh{}T'_j(i,E,{T_i}_{(2)})
\tensor{R}1\smsh{}T''_j} \qquad \\
& = &
(1\smsh{}T'_j -
\sum_{(T_{i})}{T_{i}}_{(1)}\cdot r\smsh{}T'_j(i,E,{T_i}_{(2)}))
\tensor{R} 1\smsh{}T''_j,
\end{eqnarray*}
which is clearly in
$\Ker(\pi\tensor{R}I)\subseteq\Ker(\pi\tensor{R}\pi)$.
Similarly, if node $i$ occurs in tree $T''_j$, then the corresponding
term of the coproduct applied to this element is
in $\Ker(I\tensor{R}\pi)\subseteq\Ker(\pi\tensor{R}\pi)$.
It follows
that $(\pi\tensor{R}\pi)\circ\Delta$ vanishes on $I$, since it
vanishes on
a generating set for it as an ideal in the algebra $R\smsh{k}\kLT{\D}$.

To show that $\cI(\kLT{\D})$ is compatible with the antiproduct of $\kLT{\D}$,
we will need to work with a restricted generating set of
$\cI(\kLT{\D})$.

\begin{lemma}\label{I0GenLemma}
Let $R$ be a commutative $k$-algebra, and let $\D=\Der(R)$.
Suppose that we are given a Leibnitz $\kLT{\D}$-module action on $R$.
Let $J$ be the ideal in $R\smsh{k}\kLT{\D}$ generated by
all elements of the form~(\ref{I0Gen})
where $T$ ranges over all labeled ordered trees whose root has only one
child.
Then $J=\cI(\kLT{\D}$.
\end{lemma}

\pf It is immediate that $J\subseteq\cI(\kLT{\D}$.
Let $T$ be any labeled ordered tree, and let $i$ be a node of $T$
which is labeled with $rE$.
We will prove that
\begin{equation}\label{IGenExp}
1\smsh{}T - \sum_{(T_i)}{T_i}_{(1)}\cdot r\smsh{}T(i,E,{T_i}_{(2)})
\end{equation}
is in $J$ by induction on the number of children of the root of
$T$.

If the root of $T$ has one child then the
element~(\ref{IGenExp}) is in $J$ by definition.

Suppose that the root of $T$ has $n+1$ children.
Let $T_0$ be the tree consisting of the tree whose root has one child,
which is
the first child of the root of the tree $T$,
in the order in which they occur in $T$.
Let $T_1$ be the tree whose root has as children
all of the other children of the root of the tree $T$ and their
descendents,
in the order in which they occur in $T$.
The root of the tree $T_1$ has $n$ children.
To show that element~(\ref{IGenExp}) is in $J$, we consider two
cases.
\begin{description}
\item[Node $i\in T_0$]\mbox{ }
\par\noindent
In this case let $T_0\cdot T_1 = T + \sum_j U_j$, where
each $U_j$ is a tree whose root has only $n$ children.
Then, since $T=T_0\cdot T_1 - \sum_j U_j$,
\begin{eqnarray*}
\lefteqn{
1\smsh{}T - \sum_{(T_i)}{T_i}_{(1)}\cdot r\smsh{}T(i,E,{T_i}_{(2)}) }
   \qquad & & \\
& = & (1\smsh{}T_0 -
  \sum_{(T_i)}{T_i}_{(1)}\cdot
r\smsh{}T_0(i,E,{T_i}_{(2)}))\cdot(1\smsh{}T_1) \\
&   & \qquad {} - \sum_j   (1\smsh{} U_j -
  \sum_{(T_i)}{T_i}_{(1)}\cdot r\smsh{}U_j(i,E,{T_i}_{(2)})).
\end{eqnarray*}
Note that, since the tree $T_0$ precedes $T_1$ in the product,
the labeled ordered subtree $T_i$ rooted at node $i$ is the same in
$T_0$ as in $U_j$.
The term
\[
(1\smsh{}T_0 -
  \sum_{(T_i)}{T_i}_{(1)}\cdot r\smsh{}T_0(i,E,{T_i}_{(2)}))\cdot
(1\smsh{}T_1)
\]
is in $J$ since the root of $T_0$ has one child and $J$ is an
ideal.
The terms
\[
1\smsh{} U_j -
  \sum_{(T_i)}{T_i}_{(1)}\cdot r\smsh{}U_j(i,E,{T_i}_{(2)})
\]
are in $J$ by induction on the number of children of the root of the
tree.

\item[Node $i\in T_1$]\mbox{ }
\par\noindent
In this case first note that
\begin{equation}\label{E5.0}
1\smsh{}T_1 - \sum_{(V)} V_{(1)}\cdot r\smsh{}T_1(i,E,V_{(2)}),
\end{equation}
where $V$ is the labeled subtree of $T_1$ rooted at node $i$,
is in $J$ by induction, since the root of $T_1$ has only $n$
children.
Let
\[
T_0\cdot T_1=T+\sum_k U_k+\sum_\ell U'_\ell,
\]
where the $U_k$ are the
trees in the product in which $T_0$ is not attached to the node $i$ or
to any of its descendents, and the $U'_\ell$ are the
trees in the product in which $T_0$ is attached to the node $i$ or to
one of its descendents.
Now the element
\begin{eqnarray}
\lefteqn{
(1\smsh{}T_0)\cdot(1\smsh{}T_1 - \sum_{(V)} V_{(1)}\cdot r
      \smsh{}T_1(i,V_{(2)},E))} \quad & & \nonumber \\
& = & 1\smsh{}T_0\cdot T_1 - \sum_{(V)} T_0\cdot V_{(1)}\cdot r\smsh{}
      T_1(i,V_{(2)},E) \nonumber \\
&   & \qquad {}
- \sum_{(V)} V_{(1)}\cdot r\smsh{} T_0\cdot T_1(i,V_{(2)},E)
      \nonumber \\
& = & 1\smsh{}T + \sum_k 1\smsh{}U_k + \sum_\ell 1\smsh{}U'_\ell
      \nonumber \\
&   & \qquad {} - \sum_{(V)}
      T_0\cdot V_{(1)}\cdot r\smsh{}T_1(i,V_{(2)},E) \nonumber \\
&   & \qquad {} - \sum_{(V)} V_{(1)}\cdot r \smsh{}T(i,V_{(2)},E)
      \nonumber \\
&   & \qquad {} - \sum_{(V)} V_{(1)}\cdot r\smsh{}U_k(i,V_{(2)},E)
      \nonumber \\
&   & \qquad {} - \sum_{(V)} {V}_{(1)}\cdot r\smsh{}T_1(i,T_0\cdot V_{(2)},E)
      \nonumber \\
& = & 1\smsh{}T - \sum_{(V)} V_{(1)}\cdot r \smsh{}T(i,V_{(2)},E)
      \label{I0Term1} \\
&   & \qquad {} + \sum_k \Big(1\smsh{}U_k - \sum_{(V)}
      V_{(1)}\cdot r\smsh{} U_k(i,V_{(2)},E)\Big) \label{I0Term2} \\
&   & \qquad {} + \sum_\ell\Big(1\smsh{}U'_\ell - \sum_{(V')}
      V'_{(1)}\cdot r\smsh{} U'_\ell(i,V'_{(2)},E)\Big)
      \label{I0Term3}
\end{eqnarray}
where $V'$ is the labeled subtree rooted at node $i$ in $U'_\ell$,
is in the ideal $J$, since the element~(\ref{E5.0}) is in the ideal.
Since the roots of the trees $U_k$ and $U'_\ell$ have only $n$ children,
by induction the terms~(\ref{I0Term2}) and~(\ref{I0Term3}) are in
$J$.
Therefore the term~(\ref{I0Term1}) is in the ideal $J$.
\end{description}

This completes the proof of the lemma.

\pfend

We now use Lemma~\ref{I0GenLemma} to show that $\cI(\kLT{\D}$
is compatible with
the antiproduct. Denote by $\pi$ the projection of $R\smsh{k}\kLT{\D}$
onto $R\smsh{k}\kLT{\D}/\cI(\kLT{\D}$.
We need to show that $(\pi\btensor\pi)\circ E(\cI(\kLT{\D})=0$.
By~\cite[Proposition 6]{NichKk}, since $R\smsh{k}\kLT{\D}$ is pointed,
it is sufficient to show that $(\pi\btensor\pi)\circ E=0$ on a
generating set of the ideal $\cI(\kLT{\D}$.
By Lemma~\ref{I0GenLemma} it is sufficient to consider the value of
$(\pi\btensor\pi)\circ E$ on elements of the form~(\ref{I0Gen}),
when the root of the tree $T$ has only one child.
In this case
\begin{eqnarray*}
\lefteqn{E(1\smsh{}T - \sum_{(T_i)}{T_i}_{(1)}\cdot r
    \smsh{}T(i,E,{T_i}_{(2)}))} \qquad \\
 & = & 1\smsh{}T\btensor1\smsh{}1 - 1\smsh{}1\btensor1\smsh{}T \\
 &   & \qquad {} - \sum_{(T_i)}{T_i}_{(1)}\cdot r
       \smsh{}T(i,E,{T_i}_{(2)})\btensor1\smsh{}1 \\
 &   & \qquad {} + \sum_{(T_i)}{T_i}_{(1)}\cdot r
       \smsh{}1\btensor1\smsh{}T(i,E,{T_i}_{(2)}) \\
 & = & (1\smsh{}T
       - \sum_{(T_i)}{T_i}_{(1)}\cdot r\smsh{}T(i,E,{T_i}_{(2)}))
       \btensor1\smsh{}1 \\
 &   & \qquad {} - 1\smsh{}1\btensor (1\smsh{}T - \sum_{(T_i)}
       {T_i}_{(1)}\cdot r\smsh{}T(i,E,{T_i}_{(2)})),
\end{eqnarray*}
and this is clearly annihilated by $\pi\btensor\pi$.
This proves Proposition~\ref{I0HopfIdlProp}.

\begin{dfn}\label{LabelBasisDef}
Suppose that the $\kLT{\D}$-module algebra structure of $R$ is Leibnitz.
Let $\B$ be an $R$-basis for $\D$.
Repeated applications of the substitution
\begin{equation}\label{SubstRule}
1\smsh{}T = \sum_u \sum_{(T_i)} {T_i}_{(1)}\cdot r_u
               \smsh{}T(i,X_u,{T_i}_{(2)}),
\end{equation}
where the node $i$ in the tree $T$ is labeled with
$\sum_u r_uX_u$, with $r_u\in R$, $X_u\in\B$,
gives a map
$\alpha_\B : R\smsh{k}\kLT{\D} \rightarrow R\smsh{k}\kLT{\B}$.
\end{dfn}

\begin{lemma}
Suppose that the $\kLT{\D}$-module algebra structure of $R$ is Leibnitz,
and let
$\B$ be an $R$-basis for $\D$.
Then the map $\alpha_\B : R\smsh{k}\kLT{\D} \rightarrow
R\smsh{k}\kLT{\B}$ given in Definition~\ref{LabelBasisDef} is
well-defined.
\end{lemma}

\pf It is sufficient to show that application of the
substitution rule~(\ref{SubstRule})
to two nodes does not depend on the order in which the rule is
applied to the nodes.
If neither node is an ancestor of the other, then it follows
immediately that the result is independent of the order in which the
rule is applied.
If one node is an ancestor of the other, then the independence of the
result of the order of application follows immediately from the fact
that the algebra module structure is Leibnitz.
This completes the proof of the lemma.

\pfend

Let $\beta_\B : R\smsh{k}\kLT{\B} \rightarrow R\smsh{k}\kLT{\D}$ be the
inclusion map.
Then $\alpha_\B\circ\beta_\B$ is the identity on $R\smsh{k}\kLT{\B}$.

\begin{lemma}\label{alphalemma}
Suppose that the $\kLT{\D}$-module algebra structure of $R$ is Leibnitz,
and let $\B$ be an $R$-basis for $\D$.
Then
\[
\Ker\alpha_\B = \cI(\kLT{\D}.
\]
\end{lemma}

\pf
Note that $\Ker\alpha_\B$ is the linear span of elements of the
form
\[
1\smsh{}T - \sum_{u,(T_i)} T_{i(1)}\cdot s_u\smsh{}T(i,X_u,T_{i(2)})
\]
where the node $i$ in the tree $T$ is labeled with
$\sum_u s_uX_u$, with $s_u\in R$, $X_u\in\B$.
These elements are all of the form~(\ref{I0Gen}), which span
$\cI(\kLT{\D}$, so that $\Ker\alpha_\B\subseteq\cI(\kLT{\D}$.

The ideal $\cI(\kLT{\D})$ is the linear span of elements of the form
\begin{eqnarray*}
\lefteqn{1\smsh{}T - \sum_{(T_i)}{T_i}_{(1)}\cdot r
   \smsh{}T(i,E,{T_i}_{(2)})} \qquad \\
 & = &
1\smsh{}T - \sum_{u,(T_i)}({T_i}_{(1)}\cdot r)({T_i}_{(2)}\cdot s_u)
   \smsh{}T(i,X_u,{T_i}_{(3)}) \\
 & = &
1\smsh{}T - \sum_{u,(T_i)}{T_i}_{(1)}\cdot (rs_u)
   \smsh{}T(i,X_u,{T_i}_{(2)}),
\end{eqnarray*}
Where $rE$, $r$, $s_u\in R$, $E\in\D$, is a factorization of the label
of a node $i$, $X_u\in\B$, and $E=\sum_u s_uX_u$, so that
$\cI(\kLT{\D}\subseteq\Ker\alpha_\B$.

This completes the proof of the lemma.
\pfend

We have the following theorem, which says that we can use
$R\smsh{k}\kLT{\B}$ to do computations in $R\smsh{k}\kLT{\D}$.

\begin{thm}\label{QuotI0Thm}
Let $R$ be a commutative $k$-algebra.
Assume that $\D=\Der(R)$ is free as an $R$-module.
Suppose that we have a Leibnitz $\kLT{\D}$-module algebra structure on
$R$. Let $\B$ be an $R$-basis of $\D$.
Then
\[
\alpha_\B:R\smsh{k}\kLT{\D}/\cI(\kLT{\D}) \cong R\smsh{k}\kLT{\B}.
\]
\end{thm}

\pf
>From Lemma~\ref{alphalemma} $\Ker\alpha_\B = \cI(\kLT{\D})$
so that the map $\alpha_\B$ in injective.

>From the fact that
$\alpha_\B\circ\beta_\B$ is the identity on $R\smsh{k}\kLT{\B}$
it follows the $\alpha_\B$ is surjective.
\pfend

If $\B$ is an $R$-basis of $\D=\Der(R)$, then there is a
bijection between the Leibnitz $\kLT{\D}$-module algebra structures of
$R$ and the $\kLT{\B}$-module algebra structures of $R$.
According to~\cite[Theorem 5.1]{GLtrees}, $\kLT{\B}$ is freely generated
as an
associative algebra by the set $\mathcal{X}$ of trees whose root has a
single child, which are labeled with elements of $\B$.
Therefore, there is a bijection between Leibnitz
$\kLT{\D}$-module algebra structures on $R$ and functions from
$\mathcal{X}$ to $\Der(R)$.

In particular, there are Leibnitz module algebra structures for the
Hopf
algebra of labeled ordered trees, labeled with elements of $\Der(R)$, on
$R=k[X_1$, \ldots, $X_N]$ other than the example given in Example~\ref{ex3}.
In particular, there exist module algebra structures for the Hopf
algebra of trees
labeled with elements of $\{\ptl{X_1}$, \ldots, $\ptl{X_N}\}$ under
which trees with more than two nodes whose root has only one child
act as non zero first-order differential operators.
We will see in the next section that under an additional hypothesis, the
$\kLT{\D}$-module algebra structure is determined by the actions of
trees with
two nodes, which correspond to the actions of the elements of $\D$, and
trees with three nodes whose roots have only a single child, which
correspond to the connection on $R$.

\section{Coherent actions and connections}\label{ConnectionSect}

In this section we discuss how certain actions of $\kLT{\D}$ on $R$ are
determined by the action of $E\in\D=\Der(R)$, and of the action of the
connection $\Conn{E}{F}$ for $E$, $F \in\D$.
Throughout this section $R$ is a commutative $k$-algebra.

We consider actions of $\kLT{\D}$ on $R$ under which $\vv(E)$ acts as
$E$, and under which $\uu(F;\vv(E))$ acts as $\Conn{E}{F}$.

\begin{dfn}\label{CohDef}
Suppose that $U$ is a labeled ordered tree whose root has a
single child and
which acts on $R$ as the differential operator $E_U$, and suppose that
$T$ is a labeled ordered tree which contains $U$ as a subtree.
Denote by $T(U|\vv(E_U))$ the labeled ordered tree resulting from
replacing the subtree $U$ with the tree $\vv(E_U)$.
The action of $\kLT{\D}$ on $R$ is called {\em coherent\/} if for all
labeled ordered trees $U$ whose root has a single child, and all
labeled ordered trees $T$ which contain $U$ as a subtree, the actions on
$R$ of the trees $T$ and $T(U|\vv(E_U))$ are identical, that is,
$T\sim T(U|\vv(E_U))$.
\end{dfn}

The actions defined in Example~\ref{ex3} and in Construction~\ref{ex4}
are coherent.

Note that $\kLT{\D}$ is isomorphic as an algebra to the free associative
algebra generated by the labeled ordered trees whose roots have only one
child (\cite{GLtrees}[Th.~5.1] --- what is called $\LOT$
in~\cite{GLtrees} is
called $\T$ here) so non coherent actions of $\kLT{\D}$ on
$R$ can be easily constructed.

\begin{thm}\label{CohDetThm}
Let $R$ be a commutative $k$-algebra, and let $\D=\Der(R)$.
Suppose a coherent action of $\kLT{\D}$ on $R$ is given.
Then the action of $\kLT{\D}$ on $R$ is completely determined by the
action $E$ of the trees $\vv(E)$, and the action $\Conn{E}{F}$ of the
trees $\uu(F;\vv(E))$ for all $E$, $F\in\D$.
\end{thm}

\pf
The proof uses two lemmas.

\begin{lemma}\label{BrushEquation}
Let $E_1$, \ldots, $E_n\in\D$.  Then
\begin{eqnarray*}
\lefteqn{\uu(F;\vv(E_1),\ldots,\vv(E_n)) = {}} \\
\quad &  & \vv(E_1)\cdot\uu(F;\vv(E_2),\ldots,\vv(E_n)) \\
      &  & \qquad {} - \sum_{i=2}^n
           \uu(F;\vv(E_2),\ldots,\uu(E_i;\vv(E_1)),\ldots,\vv(E_n)) \\
      &  & \qquad {}-\ttt(\vv(E_1),\uu(F;\vv(E_2),\ldots,\vv(E_n)).
\end{eqnarray*}
\end{lemma}

\pf The proof of this lemma follows immediately from the
definition of multiplication for trees.
\pfend

\begin{lemma}\label{BrushAction}
Suppose that $\kLT{\D}$ acts coherently on $R$.
Let the action of the tree $\uu(F;\vv(E_2),\ldots,\vv(E_n))$ on $R$ be
denoted by $G\in\D$, and
let the action of the tree $\uu(E_i;\vv(E_1))$ on $R$ be denoted by
$H_i\in\D$.  Then
\begin{eqnarray*}
\lefteqn{\uu(F;\vv(E_1),\ldots,\vv(E_n)) \sim {}}\\
\qquad & & \vv(E_1)\cdot\vv(G) \\
       & & {} -
           \sum_{i=2}^n \uu(F;\vv(E_2),\ldots,\vv(H_i),\ldots,\vv(E_n))
           - \ttt(\vv(E_1),\vv(G)).
\end{eqnarray*}
\end{lemma}

\pf This lemma follows immediately from Lemma~\ref{BrushEquation} and
from the definition of coherence.
\pfend

We now prove Theorem~\ref{CohDetThm}.
The action of labeled trees with two nodes is determined by the action
of
$\D$ on $R$.
Repeated application of the definition of coherence shows
that the
action of trees with more than two nodes in $\T(\D)$ is determined by
the action of trees of the form
$\uu(F;\vv(E_1),\ldots,\vv(E_n))$.
We prove by induction on $n$ that this action is determined by the
actions of $E$ and $\Conn{F}{E}$ for all $E$, $F\in\D$.

For $n=1$ this is simply the assertion that the action of
$\uu(E;\vv(F))\sim\Conn{F}{E}$ is determined.
Suppose that the action is determined for $n$.
We prove that it is determined for $n+1$.
Lemma~\ref{BrushAction} implies that the action of a tree of the form
$\uu(F;\vv(E_1),\ldots,\vv(E_{n+1}))$ is determined by the action of
trees of the form $\vv(E)$, of trees of the form
$\uu(F;\vv(E_1),\ldots,\vv(E_n))$, and of trees of the form
$\ttt(\vv(E),\vv(F))$.
The action of trees of the first form is given by hypothesis.
The action of trees of the second form is determined by the induction
hypothesis.
The action of trees of the third form is determined by
Lemma~\ref{PairEquation}. This completes the proof of the theorem.

\pfend

\end{document}